\documentclass[10pt,twoside]{amsart}
\usepackage{amssymb,graphicx,verbatim}
\newtheorem{thm}{Theorem}
\newtheorem{lem}[thm]{Lemma}
\newtheorem{prop}[thm]{Proposition}

\newtheorem{defi}[thm]{Definition}
\newcommand{\R}{\mathbb{R}}

\newcommand{\C}{\mathbb{C}}

\newcommand{\finprf}{\unskip\null\hfill$\square$\vskip 0.3cm}

\newcommand{\be}{\begin{equation}}
\newcommand{\ee}{\end{equation}}

\begin{document}


\title[Self-adjointness via Hardy-like inequalities]{Self-adjointness via partial Hardy-like inequalities}

\author[M.J. Esteban]{Maria J. Esteban$^1$}
\address{$^1$Ceremade, Universit\'e Paris Dauphine, Place de Lattre de Tassigny, F-75775 Paris C\'edex 16, France}
\email{esteban@ceremade.dauphine.fr}

\author[M. Loss]{Michael Loss$^2$}
\address{$^2$School of Mathematics, Georgia Institute of Technology, Atlanta, GA 30332, USA}
\email{loss@math.gatech.edu}

\begin{abstract}
Distinguished selfadjoint extensions of operators which are not semibounded can be deduced from the positivity of the Schur Complement
(as a quadratic form). In practical applications this amounts to proving a Hardy-like inequality. Particular cases are Dirac-Coulomb operators where distinguished selfadjoint extensions are obtained for the optimal range of coupling constants.\end{abstract}
\subjclass[2000]{}
\date{\today}

\maketitle
\vspace{2pc}\noindent{\it Keywords}: Relativistic quantum mechanics, Dirac operator, self-adjoint operator, self-adjoint extension, Schur complement.

\thispagestyle{empty}


\section{Introduction.}

In \cite{Esteban-Loss-07} we defined distinguished self-adjoint extensions of Dirac-Coulomb operators in the optimal range  for the coupling constant. This was done by using a Hardy-like inequality which allowed the extension of one component of the operator by using the Friedrichs extension. Then, the remaining component could be extended by choosing the right domain for the whole operator. The method of proof used simple arguments of distributional differentiation. This work was the sequel of a series of papers where distinguished self-adjoint extensions of Dirac-Coulomb like operators were defined by different methods almost in the optimal range, without reaching the limit case (see 
\cite{Schmincke-72, Schmincke-72B, Wust-73, Wust-75, Wust-77, Nenciu-76, Klaus-Wust-78}).

Here we present an abstract version of the method introduced in \cite{Esteban-Loss-07}. We believe that this will clarify the precise structure and hypotheses necessary to define distinguished self-adjoint extensions by this method.

The main idea in our method is that Hardy-like inequalities are fundamental to define distinguished (physically relevant) self-adjoint extensions {\it even}  for operators that are not bounded below. 

We are going to apply our  method to operators $H$ defined on ${\mathcal D}_0^2$, where $\,{\mathcal D}_0\,$ is some dense subspace of a Hilbert space $\,{\mathcal H_0}$. The general structure taken into account here is:
\be\label{operator} H= \left(\begin{matrix}P & Q \\ T & -S \end{matrix}\right)\,,
\ee
where all the above operators satisfy    $Q=T^*$, $P=P^*, \,S=S^*$ and  $S\ge c_1 I >0$. Moreover we assume that $\,P, \,Q, \,S, \,,T, \,S^{-1}T\,$ and   $\,QS^{-1}T\,$ send 
$\,{\mathcal D}_0\,$ into  $\,{\mathcal H_0}  $. 

In the Dirac-Coulomb case our choice was $\, {\mathcal H}_0= L^2(\R^3, \C^2)\,$ and 
$$ P=V+2-\gamma,\; Q=T= -i\sigma\cdot\nabla, \; S = \gamma-V\,,$$
where $V$ is a potential bounded from above satisfying 
\begin{equation}\label{sotto2}
\,\sup_{x\ne 0}\, |x||V(x)| \le 1\,.
\end{equation}
Moreover, $\sigma_i,\, i=1,2,3, $ are the Pauli matrices (see \cite{Esteban-Loss-07}) and $\,\gamma\,$ is a constant slightly above $\,\max_{\R^3}V(x)$. For $\,{\mathcal D}_0\,$ we chose 
$C^\infty_c(\R^3,\C^2)  $.
Note that in our paper \cite{Esteban-Loss-07}, where we deal with Dirac-Coulomb like operators, there is an omission. We forgot to specify the conditions on the potential $V$ so that 
$QS^{-1}T$ is a symmetric operator on $C^\infty_c(\R^3,\C^2)$.
The natural condition is that each component of
\begin{equation}
(\gamma - V)^{-2} \nabla V
\end{equation}
is locally square integrable. This is easily seen to be true for the Coulomb-type potentials.

\medskip
In the general context of the operator $H$, as defined in \eqref{operator}, our main assumption is that there exists a constant $\,c_2 >0\,$ such that for all $\,u\in {\mathcal D}_0$,
\be\label{R7}
q_{c_2}(u,u) :=\big((S+c_2)^{-1}Tu, Tu\big)+\big((P-c_2)u,u\big)\geq 0\,.
\ee
Note that since $\, \frac{d}{d\alpha}q_\alpha(u,u)\leq -(u,u)$, \eqref{R7} implies in fact that  for all $\, 0\le \alpha \le c_2\,$ and for all $u\in {\mathcal D}_0$,
\be\label{R75}
q_\alpha(u,u) :=\big((S+\alpha)^{-1}Tu, Tu\big)+\big((P-\alpha)u,u\big)\geq 0\,.
\ee

\medskip
Antoher consequence of assumption \eqref{R7} is that the quadratic form
\be\label{beegamma}
q_0(u,u) = \big(S^{-1}Tu, Tu\big)+\big(Pu,u\big)\,,
\ee
defined for $\,u\in {\mathcal D}_0$, is positive definite:
\be
q_0(u,u) = \big(S^{-1}Tu, Tu\big)+\big(Pu,u\big) \geq c_2(u,u).
\ee
  Note that the operator $\,P+QS^{-1}T$ which is associated with the quadratic form $q_0$ is actually the Schur complement of $-S$. {Note also that by our assumptions on $P,\, Q, \,T\, S$ and by
\eqref{R7}, for any $0 \le \alpha \le c_2$, $q_\alpha$ is the quadratic form associated with a positive symmetric operator. Therefore, by Thm. X.23 in \cite{Reed-Simon-78}),} it is closable and we denote its closure by $\widehat{q}_\alpha$ and its form domain,
which is easily seen to be independent of $\alpha$
(see \cite{Esteban-Loss-07})  by
$\mathcal{H}_{+1}$.  
Our main result states the following:

\begin{thm}\label{thmain}
Assume the above hypotheses on the operators $P, Q, T, S$ and \eqref{R7}. Then there is a unique self-adjoint extension of $H$ such that the domain of the operator is contained in $\mathcal{H}_{+1}\times {\mathcal H_0}$.
\end{thm}

\noindent{\bf Remark. } Note that what this theorem says that ``in some sense" the Schur complement of $-S$ is positive, and therefore 
has a natural self-adjoint extension, then one can define a distinguished self-adjoint extension of the operator $H$ which is unique among those whose domain is contained in the form domain of the Schur complement of $-S$ times ${\mathcal H_0}$. 

\section{Intermediate  results and proofs.}

 We denote by $R$ the unique selfadjoint operator associated with $\widehat{q}_0$: for
all $u \in D(R) \subset \mathcal{H}_{+1}$,
\be
\widehat{q}_0(u, u) = (u, R u) \ .
\ee
$R$ is an isometric isomorphism
from ${\mathcal H}_{+1}$ to its dual ${\mathcal H}_{-1}$. Using the second representation theorem in \cite{Kato}, Theorem
2.23, we know that
${\mathcal H}_{+1}$ is the operator domain of $R^{1/2}$,
and
\be
\widehat{q}_0(u, u) = (R^{1/2}u, R^{1/2} u) \ ,
\ee
for all $u \in  {\mathcal H}_{+1}$.

\begin{defi} \label{defi-domain} We define
the domain $\mathcal{D}$ of $\,H\,$ as the collection of all pairs
$u \in {\mathcal H}_{+1}$, $ v \in {\mathcal H_0}\,$ such that
\be \label{distributional}
Pu+Qv \ , \quad   Tu-Sv  \in \ {\mathcal H_0} \ . 
\ee
 
\end{defi}
The meaning of these two expressions is in the weak (distributional) sense, i.e., the linear functional $\,(P\eta, u) +(Q^*\eta, v)\,$, which is defined for all test functions $\eta \in \mathcal{D}_0$, extends uniquely to a bounded linear functional on $\,{\mathcal H_0}$. Likewise the same for 
$\,(-S\eta, v) +(T^*\eta, v)\,$. 

On the domain $\mathcal{D}$, we define the
 operator $H$ as
\be
H{u\choose v}= {Pu+Qv \choose Tu-Sv}\,.
\ee
Note that for all vectors $(u,v) \in \mathcal{D}$ the
expected total energy is finite.

The following two results are important in the proof of Theorem \ref{thmain}.

\begin{prop}\label{prop-scale}
Under the assumptions of Theorem \ref{thmain}
 \be\label{domain-estimate}
{\mathcal H}_{+1}\subset\Big\{u\in {\mathcal H_0}\,:\; S^{-1}Tu\in {\mathcal H_0}\Big\}\,,
\ee
 where the embedding holds in the continuous sense.
Therefore, we have the `scale of spaces' \; ${\mathcal H}_{+1} \subset
{\mathcal H_0} \subset {\mathcal H}_{-1}$.
 \end{prop}
 
 \proof
Choose $\,c_2 \ge \alpha>0\,$. Since $\,S \ge c_1 I\,$, we have for all $\,0<\delta\leq \frac{c_1 \alpha}{c_1+\alpha}\,$ 
\be
S^{-1}-(S+\alpha)^{-1}\geq \delta\,S^{-2}\,,
\ee
and so, for all $u\in {\mathcal D}_0$,
\be 
q_0 (u,u)\geq q_\alpha (u,u) + \alpha \,(u,u) + \delta\,(S^{-1}Tu, S^{-1}Tu) \ge \delta \,(u,u) + \delta\,(S^{-1}Tu, S^{-1}Tu) \ .
\ee 
The proof can be finished by density arguments.
\finprf

\begin{lem} \label{lem-noidea}
For any $F$ in ${\mathcal H_0}$, 
\be
QS^{-1} F \in \mathcal{H}_{-1} \ .
\ee
\end{lem}

\proof
By our assumptions on $H$ and by Proposition \ref{prop-scale}, for every
$\eta \in \mathcal{D}_0$,
\be
\left|(S^{-1}T\eta, \, F)\right|
 \le \delta^{-1/2}\,\Vert \eta \Vert_{\mathcal{H}_{+1}} \, \Vert F\Vert_2 \ .
\ee
Hence, the linear functional
\be
\eta \to (Q^*\eta, S^{-1}F)
\ee
extends uniquely to a bounded linear functional on $\mathcal{H}_{+1}$.
\finprf

\noindent{\sl Proof of Theorem \ref{thmain}.}
We shall prove Theorem \ref{thmain} by showing that $\,H\,$ is symmetric and  a bijection
from its domain $\mathcal{D}$ onto ${\mathcal H_0}$. 
To prove the symmetry
we have to show that for both pairs $(u,  v)$, $(\tilde{u}, \tilde{ v})$ in the domain ${\mathcal D}$,
\be
\left(H\,{u \choose v}, \, {\tilde u \choose \tilde v}\right)
= (Pu+Qv, \,\tilde u) + (Tu-Sv, \,\tilde v)
\ee
equals
\be
(u, \,P\tilde u+Q\tilde v) + (v,\, T\tilde u-S\tilde v)= \left({u \choose v}, \, H\,{\tilde u \choose \tilde v}\right)
.
\ee
First, note that since $(u, v)$ is in the domain,
\be \label{elltwo}
S(v-S^{-1}Tu) \in  {\mathcal H_0}) \ .
\ee 
We now claim that
\be
(Pu+Qv, \,\tilde u)
=(R u, \tilde{u}) +(\,S(v-S^{-1}Tu), \,S^{-1}T\tilde u\,)\,.
\ee
Note that each term makes sense. The one on the left, by definition
of the domain and the first on the right, because both $u, \tilde{u}$ are in ${\mathcal H}_{+1}$. The second term on the right
side makes sense because of \eqref{elltwo} above and Proposition \ref{prop-scale}.
Moreover both sides coincide for $\tilde{u}$ chosen to be a test function and both are continuous in $\tilde{u}$ with 
respect to the ${\mathcal H}_{+1}$ -norm. Hence the two expressions coincide on the domain. Thus we get
that
\be
\left(H\,{u \choose v}, \, {\tilde u \choose \tilde v}\right)
\ee
equals
\be
(R u, \tilde{u}) - (\,S(v-S^{-1}Tu), \, \tilde v-S^{-1}T\tilde u \,)\,,
\ee
an expression which is symmetric in $(u, v)$ and $(\tilde{u}, \tilde{ v})$.
To show that the operator is onto,  pick any $F_1, F_2$ in ${\mathcal H_0}$. Since $R$ is an isomorphism, there exists
a unique
$u$ in ${\mathcal H}_{+1}$ such that
\be \label{phiequation}
R u = F_1 +QS^{-1}F_2\ .
\ee
Indeed, $F_1$ is in ${\mathcal H_0}$ and therefore in ${\mathcal H}_{-1}$. Moreover the second term is also in ${\mathcal H}_{-1}$ by Lemma \ref{lem-noidea}. 

Now define $ v$ by
\be \label{chiequation}
 v = S^{-1}(Tu-F_2)\,,
\ee
which by Proposition \ref{prop-scale} is in ${\mathcal H_0}$. 

Now for any test function $\eta$ we have that
\be
(P\eta, \, u) + (Q^*\eta, \, v) =(P\eta,\, u) + (T\eta, \, v) =
(P\eta,\, u) +(\,T\eta,\,S^{-1}Tu)+(T\eta, \, (v-S^{-1}Tu)\,)
\ee
which equals
\be
(\eta,Ru) + (T\eta, \,( v -S^{-1} T u)
= (\eta, F_1)
\ee
This holds for all test functions $\eta$, but since $F_1$ is in ${\mathcal H_0}$,
the functional
$\eta \to (P\eta, \,u) + (T\eta, \, v)$
extends uniquely to a linear continuous functional on ${\mathcal H_0}$ which
implies that
\be
Pu+Qv = F_1 \ .
\ee
Hence $(u,  v)$ is in the domain $\mathcal{D}$ and the operator $\,H\,$ applied to $(u,  v)$ yields $(F_1, F_2)$.

Let us now prove the injectivity of $\,H$. Assuming that
\be\label{esc3}
H{u\choose v}={0 \choose 0}\,,
\ee
we find by \eqref {phiequation} and \eqref {chiequation}, 
$$ v=S^{-1}Tu\;,\qquad R u=0\,.
$$
Since $R$ is an isomorphism, this implies that $u= v=0$.

It remains to show the uniqueness part in our theorem. By the bijectivity result proved above, for all ${F_1\choose F_2}\in  {\mathcal H_0}^2$, there exists a unique pair $\,(\hat u \,, \hat v)\in \mathcal{H}_{+1}\times {\mathcal H_0}\,$ such that $\,H{u\choose v} = {F_1\choose F_2}$. Let us now pick any other self-adjoint extension with domain $\,{\mathcal D}'$ included in $\,\mathcal{H}_{+1}\times {\mathcal H_0}\,$. Then  for all $\,(u,v)\in {\mathcal D}'$, $\,H{u \choose v}\,$ belongs to ${\mathcal H_0}^2$. Hence there exist a unique pair $\,(\hat u \,, \hat v)\in \mathcal{H}_{+1}\times {\mathcal H_0}\,$ such that $\,H{\hat u\choose \hat v} =H{u\choose v}$. But, by the above considerations on injectivity, $u=\hat u$  and $v= \hat v$. Therefore, ${\mathcal D}'\subset {\mathcal D}\,$ and so necessarily, ${\mathcal D}'= {\mathcal D}$.
 \finprf

\bigskip\noindent{\bf Acknowledgments.} M.J.E. would like to thank M. Lewin,  E. S\'er\'e and J.-P. Solovej for various discussions on the self-adjointness of Dirac operators.

\smallskip
M.J.E. and M.L. wish to express their gratitute to  Georgia Tech and Ceremade for their hospitality. M.J.E. acknowledges support from ANR Accquarel project and European Program ``Analysis and Quantum'' HPRN-CT \# 2002-00277. M.L. is partially supported by U.S. National Science Foundation grant DMS DMS 06-00037.

\medskip\noindent\copyright\, 2007 by the authors. This paper may be reproduced, in its entirety, for non-commercial purposes.


\end{document}